\documentclass{cmsams-l}

\numberwithin{equation}{section}

\def\R{{\car}}

\def\Re{{\br}}

\def\Ze{{\bz}}

\def\Co{{\bc}}

\def\H{{\ch}}
\def\<{\leq_\Re}
\def\sse{{\it\thinspace iff\thinspace\ }}
\def\ad{{\rm ad}\thinspace}

\def\D{\Delta}\def\d{\delta}\def\l{\lambda}
\def\s{\sigma}\def\t{\vartheta}
\def\e{\varepsilon}\def\b{\beta}

\def\cd{{\mathcal D}}

\def\ch{{\mathcal H}}

\def\cam{{\mathcal M}}
\def\cn{{\mathcal N}}

\def\car{{\mathcal R}}

\def\cz{{\mathcal Z}}

\def\bc{{\mathbb C}}

\def\br{{\mathbb R}}
\def\bz{{\mathbb Z}}
\newtheorem{Thm}{Theorem}[section]

\newtheorem{Prop}[Thm]{Proposition}
\newtheorem{Lemma}[Thm]{Lemma}
\theoremstyle{definition}
\newtheorem{Dfn}[Thm]{Definition}

\theoremstyle{remark}
\newtheorem{rem}[Thm]{Remark} 
\newtheorem{ack}{Acknowledgement} 
\begin{document}
\title{Inclusions of second quantization algebras}
\author{Franca Figliolini} 
\author{Daniele Guido}
\address{D. Guido, Dip. di Matematica, 
Universit\`a della Basilicata, I--10085 Potenza, Italy.}
\email{guido@unibas.it}
\thanks{The second named author is supported in part by MURST and EU}
\maketitle

\section{Introduction.}
In this note we study inclusions of second quantization algebras,
namely inclusions $\R(M_0)\subset\R(M_1)$ of von~Neumann algebras on
the Fock space $e^\H$ ($\H$ is a separable complex Hilbert space)
generated by the Weyl unitaries with test functions in the closed,
real linear subspaces $M_0$, $M_1$ of $\H$.  More precisely we concentrate
our attention to the case where both $\R(M_0)$ and $\R(M_1)$ are
standard w.r.t. the vacuum vector $e^0\in e^\H$, since in this case
the tower and tunnel associated with the inclusion (and the
corresponding relative commutants) can themselves be realized as second
quantization algebras on the same space:
$$
\dots\subset\R(M_{-1})\subset\R(M_{0})\subset\R(M_{1})\subset\R(M_{2})\dots
$$

First we show that the class of irreducible inclusions of standard
second quantization algebras is non empty, and that they are depth two
inclusions, namely $\R(M_0)'\cap\R(M_3)$ is a factor.  Then we prove
that, when $M_0\subset M_1$ is a (not necessarily irreducible)
inclusion of standard spaces with finite codimension $n$, $\R(M)$ is
isomorphic to the cross product of $\R(N)$ with $\br^n$.  On the
contrary, when the codimension is infinite, we show that the inclusion
may be non regular (see subsection \ref{4.1}).

Second quantization algebras and their inclusions occur when studying
algebras of local observables for the free fields.  Inclusions of
local observable algebras are in general neither irreducible nor come
from a finite codimension inclusion of real vector spaces. In
\cite{GLWi1} however, local algebras for conformal field theories on
the real line are studied, and it is shown that the inclusion of the
real vector space corresponding to a bounded interval for the $n+p$-th
derivative of the current algebra into the real vector space for the
same interval and the $n$-th derivative theory has codimension $p$
(and is irreducible when $p=1$). We show in Theorem \ref{(2.1)} that
the corresponding inclusion of second quantization algebras is given
by a cross product for any $n\geq0$, $p>0$.

Our analysis was also motivated by results concerning depth two
inclusions of von~Neumann algebras. It is well known that, analyzing
the Jones' tower associated with the inclusion of a von~Neumann
algebra $\cn$ into the cross product $\cam$ of the same algebra with
an outer action of locally compact group, the family of relative
commutants presents some characteristic features, namely the first
relative commutant is trivial, the second is abelian, and the third is
a type I factor. The study of inclusions with these properties,
started in \cite{HeOc1}, was recently enriched by a pair of papers
\cite{EnNe1,Enoc1} where these properties, together with the existence
of an operator valued weight from $\cam_1$ to $\cam_0$, where shown to
characterize the inclusions given by cross products with locally
compact groups.  We show that inclusions of second quantization
algebras may produce examples of irreducible depth two inclusions with
type III third relative commutant.  By a result in \cite{EnNe1}, these
inclusions are not regular, since in that case the third relative
commutant would be type I, hence do not correspond to a
crossed--product with a locally compact group.

On the other hand, our examples satisfy many of the features of a
crossed--product, the locally compact group being replaced by an
infinite-dimensional vector space, thus furnishing examples in order
to develop a theory of non-locally compact cross products.

\section{Preliminaries}

In this section $\H$ will be a separable Hilbert space and $e^\H$
the {\it symmetric Fock space} over it, i.e. the space
$$
e^\H=\bigoplus^\infty_{n=0}\H^{\otimes_Sn}
$$
where $\H^{\otimes_Sn}$ is the subspace of the $n-$th  tensor product
of $\H$ which is pointwise invariant under the natural action of the
permutation group.

The set of {\it coherent} vectors in $e^\H$ consists of the vectors
$$e^h=\bigoplus^\infty_{n=0}\frac{h^{\otimes n}}{\sqrt{n!}}\quad.$$
This set turns out to be total in $e^\H$ (see e.g. \cite{Guic1}, p.32).

There are two important classes of operators acting on $e^\H$:
\\
{\it Second quantization operators}
$$
e^A=\bigoplus^\infty_{n=0}A^{\otimes n}\ ,
$$
where $A$ is a densely defined, closed operator on $\H$, and
\\
{\it Weyl unitaries}, which are the range of the map 
$$
h\to W(h)
$$ 
from $\H$ to the unitaries on $e^\H$ defined by

\begin{align}
W(h)e^0
&=\exp\left(-{1\over4}\|h\|^2\right)e^{{i\over\sqrt2}h},\quad h\in\H
\label{vacuum}
\\
W(h)W(k)
&=\exp\left(-{i\over2}Im(h,k)\right)W(h+k)\quad h,k\in\H
\label{CCR}
\end{align}

The vector $e^0$ is called {\it vacuum} and the relations in the
last equality  are called Canonical Commutation Relations. Via the
preceding equalities $W(h)$ becomes a well defined, isometric and
invertible (with inverse $W(-h)$) operator on the dense set spanned
by coherent vectors, and hence it extends to a unitary on $e^\H$. Weyl
unitaries generate the so-called {\it second quantization algebras}.
With each closed real linear subspace $K\subset\H$ (in the following
we shall write $K\<\H$), a von Neumann algebra $\R(K)$ is associated,
defined by 
$$
\R(K)=\{W(h),\quad h\in K\}''.
$$ 
  
In the following Proposition we recall the basic facts about these
algebras \cite{Arak1,EcOs1}: 
\begin{Prop}\label{(1.1)} Let us consider the map $K\to\R(K)$
where $K\<\H$. Then 
\item{(i)}The mentioned map is an isomorphism of
complemented nets where ``\thinspace$\wedge$'' is the intersection,
``\thinspace$\vee$'' is the generated real subspace, resp. von
Neuman algebra, and the complementation ``\thinspace$'$'' is the polar
space w.r.t. $Im(\cdot,\cdot)$, resp. the commutant. 
\item{(ii)} The map preserve the standard property, namely the vacuum
vector is standard, i.e. cyclic and separating, for $\R(K)$ if and
only if $K$ is standard, i.e.  $K{\cap}iK=0$ and $K+iK$ is dense.
Moreover the Tomita operator $S$ associated with a standard algebra
$\R(K)$ is the second quantization of the Tomita operator $s$ on $K$
defined by
$$
s:h+ik\to h-ik\qquad h,k\in K
$$
and also $J=e^j$, $\D=e^\d$, where $S=J\D^{1/2}$, $s=j\d^{1/2}$.
\end{Prop}

As a consequence of the preceding proposition, the properties of the
second quantization algebras can be studied in terms of the generating
one-particle subspaces, and in the following we shall attribute to a
real linear subspace $K$ of $\H$ all the properties of its second
quantization algebra, for instance we will say that K is a standard
type III subspace if $\R(K)$ is, etc.  We shall also say that
$M'=\{x\in\H:Im(x,y)=0\ \forall x\in M\}$ is the commutant of $M$.  In
particular, if $M_0\subset M_1$ is an inclusion of standard subspaces
of $\H$, the tower and the tunnel \cite{Jone1,PiPo1} generated by the
inclusion $\cam_0\equiv\R(M_0)\subset\cam_1\equiv\R(M_1)$ can be
spatially realized as
\begin{align*}
\cam_{i+1}=J_i \cam'_{i-1}J_i\qquad &i\geq1\cr
\cam_{i-1}=J_i \cam'_{i+1}J_i\qquad &i\leq0,
\end{align*}
where $J_i$ is the Tomita conjugation for $\{\cam_i,e^0\}$.
Therefore they are the second quantization of the tower, resp. tunnel, of
standard subspaces:

\begin{align*}
M_{i+1}=j_i M'_{i-1}\qquad &i\geq1\cr
M_{i-1}=j_i M'_{i+1}\qquad &i\leq0.
\end{align*}

\section{Inclusions of standard subspaces}

\begin{Dfn}\label{(1.1a)} We shall say that a pair $(E,F)$, $E,F\<\H$,  
is a standard pair if $E\wedge F=\{0\}$ and $E\vee F=\H$. A
standard pair will be called strongly standard if, moreover,
$E+F=\H$. If $(E,F)$ is a standard pair, we define the operator $s_{E,F}$ as
follows:
$$
\begin{matrix}
	s_{E,F}:&E+F&    \to&E+F\cr
	        &e+f&\mapsto&e-f
\end{matrix}
$$
 \end{Dfn} 

Let us observe that $N$ is standard if and only if $(N,iN)$ is a standard pair.
Also,  $N$ is a factor if and only if $(N,N')$ is a standard pair. Indeed,
if $N$ is a factor then $N\cap N'=\{0\}$, and, taking the
commutant, $N\vee N'=(N\cap N')'=\H$.

\begin{Lemma}\label{newlemma0} Let $(E,F)$ be a standard pair. Then:
\item{$(i)$}    $s_{E,F}$ is closed
\item{$(ii)$}   $s_{E,F}$ is bounded \sse $(E,F)$ is a strongly standard pair.
\item{$(iii)$} The pair $(E,F)$ can be properly extended by a standard
pair $(\tilde E,\tilde F)$ (i.e. $E\subseteq\tilde E$,
$F\subseteq\tilde F$ where at least one inclusion is proper) \sse
$(E,F)$ is not strongly standard.
\end{Lemma}

\begin{proof} 
$(i)$ Since the graph norm of $e+f$ w.r.t. $s_{E,F}$ is
$(\|e+f\|^2+\|s_{E,F}(e+f)\|^2)^{1/2}=\sqrt2(\|e\|^2+\|f\|^2)^{1/2}$, the
result follows.

$(ii)$ Immediate by the closed graph theorem.

$(iii)$ Suppose $(E,F)$ is strongly standard and let  $(\tilde E,\tilde F)$ be a
standard extension. For any $y\in\tilde E$, $y=e+f$, hence 
$f=y-e\in F\cap\tilde E\subset\tilde F\cap\tilde E$, which implies $f=0$, i.e.
$y\in E$, hence $\tilde E=E$.
 \end{proof}

\begin{Lemma}\label{newlemma1} Let $N\<\H$ be a standard subspace. 
Then
\item{$(i)$} There exists a standard  subspace $M$ such that
$N\subset M$ if and only if $(N,iN)$ is not strongly standard.  
\item{$(ii)$} If $y\in(N+iN)^c$, then $M=\{x+\l y:x\in N, \l\in \br\}$ is
standard (but not strongly standard).
\end{Lemma}

\begin{proof} If $(N,iN)$ is strongly standard and $M\supset N$, then
$(M,iM)$ cannot be standard by Lemma~\ref{newlemma0}, $(iii)$.

Now let $(N,iN)$ be not strongly standard, $y\in(N+iN)^c$ and
 $M=\{x+\l y:x\in N, \l\in \br\}$. If $z\in M\cap iM$ then
$z=n_1+\l_1 y=i n_2+i\l_2 y$ for some $n_1,n_2\in N$, $\l_1,\l_2\in\br$.
However $\l_1$ and $\l_2$ should be zero, otherwise
$y=(i\l_2-\l_1)^{-1}(n_1-in_2)\in N+iN$, hence $z\in N\cap iN=\{0\}$.
Finally, if $M$ is standard then $(N,iN)$ is not strongly standard hence
$s_{N}$ is unbounded. Since $s_N\subset s_M$, $s_M$ is unbounded, hence $M$
is not strongly standard.
\end{proof}

\begin{Lemma}\label{newlemma2} Let $N$ be a factor subspace of $\H$. Then
\item{$(i)$} There exists a subspace $M$ such that
$N\subset M$ is irreducible an irreducible 
if and only if $(N,N')$ is not strongly standard.  
\item{$(ii)$} If $y\in(N+N')^c$, and $M=\{x+\l y:x\in N, \l\in \br\}$,
then $N'\cap M=\{0\}$ (and $(M,M')$ is not strongly standard). 
\end{Lemma}

\begin{proof} If $(N,N')$ is strongly standard and $M\supset N$, then
$(M,N')$ cannot be standard by Lemma~\ref{newlemma0}, $(iii)$.

Now let $(N,N')$ be not strongly standard, $y\in(N+N')^c$ and
 $M=\{x+\l y:x\in N, \l\in \br\}$. If $z\in M\cap N'$ then
$z=n+\l y=n'$ for some $n\in N$, $n'\in N'$, $\l\in\br$.
However $\l$ should be zero, otherwise
$y=\l^{-1}(n'-n)\in N+N'$, hence $z\in N\cap N'=\{0\}$.
Finally, observe that $(M,M')$ is properly extended by $(M,N')$, hence 
$s_(M,M')$ cannot be bounded, which means that
$(M,M')$ is not strongly standard.
\end{proof}

\begin{Prop}\label{extensions} Let $N$ be a standard subspace of $\H$. Then
\item{$(i)$} There exists a standard subspace $M\supset N$ if and only if
$0\in\s(\d_N)$.
\item{$(ii)$} Suppose $N$ is a factor. There exists a subspace $M$
such that $N\subset M$ is an irreducible inclusion if and only if
$1\in\s(\d_N)$.
\item{$(iii)$} Suppose $N$ is a factor. There exists a standard
subspace $M$ such that $N\subset M$ is an irreducible inclusion if and
only if $\{0,1\}\subset\s(\d_N)$.  
\item{$(iv)$} If $N\subset M$ is an irreducible inclusion of standard
spaces,then $N$ and $M$ are type III$_1$ factors.
\end{Prop}

\begin{proof}
$(i)$ By Lemmas \ref{newlemma0}, $(ii)$, \ref{newlemma1}, $(i)$, there
exists a standard subspace $M\supset N$ \sse $s_N$, hence
$\d_N^{1/2}$, is unbounded. Since $j_N\d_N^{1/2}j_N=\d_N^{-1/2}$, this
is is equivalent to $0\in\s(\d_N^{1/2})$, which in turn corresponds to
$0\in\s(\d_N)$ by spectral mapping.

$(ii)$ Let $\Theta$ be the selfadjoint operator defined by
$\cos\Theta=\frac{|\d-I|}{\d+I}$, $\d=\d_N$, and let
$s_{N,N'}=ud^{1/2}$ be the polar decomposition of $s_{N,N'}$. By
$s_{N,N'}^2\subset I$ one obtains $u^2=1$ and $udu=d^{-1}$, and by
\cite{FiGu2}, Proposition 3.2, one gets
$\sin\Theta=\frac{|d-I|}{d+I}$.  This implies that there exists a
subspace $M$ such that $N\subset M$ is irreducible $\iff$ $s_{N,N'}$
is unbounded, $\iff$ $0\in\s(d)$, $\iff$ $\pi/2\in\s(\Theta)$, $\iff$
$1\in\s(\d)$.

$(iii)$ By the relations among $\d$, $d$ and $\Theta$ one gets that
$x\in\cd(\d^{1/2})^c$ $\iff$
$\chi_{[0,\e]}(\Theta)x\in\cd(\d^{1/2})^c$ and $y\in\cd(d^{1/2})^c$
$\iff$ $\chi_{[\pi/2-\e,\pi/2]}(\Theta)y\in\cd(d^{1/2})^c$, therefore
if $x\in\cd(\d^{1/2})^c$ and $y\in\cd(d^{1/2})^c$, then
$z=\chi_{[0,\pi/4]}(\Theta)x+\chi_{[\pi/4,\pi/2]}(\Theta)y
\in(\cd(\d^{1/2})\cup\cd(d^{1/2}))^c$, hence $M=\{x+\lambda z:x\in N,
\ \lambda\in\br\}$ has the required properties by Lemmas
\ref{newlemma1}, \ref{newlemma2}.

$(iv)$ As shown before, $N$ is a factor if and only if
$1\in\s(\d_N)$. By \cite{FiGu2}, Proposition 4.5, this implies that
$N$ is an injective type III$_1$ factor. Since $M\subset j_MN'j_M$ is
an irreducible inclusion, the same reasoning applies to $M$.
\end{proof}

\begin{rem}\label{(1.3)} If the inclusion of standard spaces
$M_0\subset M_1$ is irreducible, all the algebras of the Jones tower
(and tunnel) are  standard type III$_1$ factors. 
\end{rem}

We have proved that the set of irreducible incusions of standard factors in $\H$
is not empty. More precisely we have shown how to costruct irreducible
inclusions of standard factors with codimension 1. Iterating the process and
taking direct sums one can construct irreducible inclusions of standard
factors with any codimension. In the following we study the structure of such
inclusions.

\begin{Lemma}\label{(1.4)} If $M$ is standard in $\H$, $(M,Ker(j+I))$
is a strongly standard pair, where $j$ is the modular conjugation
of $M$.
\end{Lemma}

\begin{proof} This proof closely follows the proof of Proposition 3.2 in
\cite{FiGu2}. The spaces $M, Ker(j+I)$ are both reduced by
the spectral projections relative to $|\log\d|$, where $\d$ is the
modular operator relative to $M$, hence it is  sufficient to study
the problem on the fiber of the representation of $\H$ as direct
integral w.r.t.  $|\log\d|$. Then, as in \cite{FiGu2}, Prop.
3.2, it is sufficient to study the case in which $K$ is a standard
subspace of $\Co^2$, hence $K$ can  be seen as generated by the
vectors 
$$
y^+=\left(
\begin{matrix}
\cos{\t\over2}\cr
\sin{\t\over2}\cr
\end{matrix}\right)\qquad
y^-=\left(
\begin{matrix}
i\cos{\t\over2}\cr
-i\sin{\t\over2}\cr
\end{matrix}\right) 
$$
and $j$ is represented by the matrix $$
\left(
\begin{matrix}
0&C\cr 
C&0\cr
\end{matrix}\right)
$$
where $C$ is the complex conjugation. Therefore
$$
Ker(j+I)=\{\left(
\begin{matrix}
a\cr
-\overline a\cr
\end{matrix}\right),\quad a\in\Re\}
$$
and it is easy to see that 
$$
|Re(h,k)|\leq{\sqrt2\over2},\quad h\in M,\ k\in
Ker(j+I),\ \|h\|=\|k\|=1,
$$
which implies 
\begin{align*}
\frac{\|h+k\|^2}{\|h\|^2+\|k\|^2}=
&\frac{\|h\|^2+\|k\|^2+2Re(h,k)}{\|h\|^2+\|k\|^2}\cr
&\geq\frac{\|h\|^2+\|k\|^2-\sqrt2\|h\|\|k\|}{\|h\|^2+\|k\|^2}\cr
&\geq\frac{\sqrt2-1}{\sqrt2},
\end{align*}
namely the graph norm of $s_{M,Ker(j+I)}$ is equivalent to the Hilbert
norm, and the strongly standard property follows by Lemma
\ref{newlemma0}, $(ii)$.
\end{proof}

If $M_0\subset M_1$ are standard subspaces and $j_0$, $j_1$ are the
modular conjugations, the tower and the tunnel of standard
subspaces are inductively defined by the equations:
$$
M_{k+1}=j_kM_{k-1}',\qquad k\in\bz.
$$
We shall use the following notations:
\begin{align}
A_{k,l}&=M_k'\cap M_l\qquad k\leq l\in\Ze\label{equ(1.1)bis}
\\
B_k&=M_{k+1}\cap Ker(j_k+I)\qquad k\in\Ze\label{equ(1.1)}
\end{align}
Now we can prove some structural results about inclusions of
standard subspaces.
\begin{Prop}\label{(1.5)} Let $M_0\subset M_1$ be an inclusion of
standard subspaces. Then the following holds:

\begin{align*}
&(i)\quad M_{k+p}=M_k+\sum_{j=0}^{p-1}B_{k+j}\qquad
&k\in\Ze,\ p\geq0\cr
&(ii)\quad A_{k-1,k+p}=A_{k-1,k}+\sum_{j=0}^{p-1}B_{k+j}\qquad
&k\in\Ze,\ p\geq0\cr
&(iii)\quad B_k\quad{\rm is\ commutative}\qquad k\in\Ze\cr
&(iv)\quad \sum_{j=1}^{2p}B_{k+j}\quad{\rm is\ a\ factor}\qquad 
&k\in\Ze,\ p>0\cr
&(v)\quad\cz\left(\sum_{j=1}^{2p+1}B_{k+j}\right)= 
\{(I+\sum_{j=1}^{p}(-1)^i j_{k+2i})x:x\in B_{k+1}\}\qquad
&k\in\Ze,\ p>0,
\end{align*}
where $\cz(P)=P\cap P'$. In particular, if the inclusion $M_0\subset
M_1$ is irreducible,
\item{$(vi)$} 
\quad $A_{k-1,k+1}=B_k$ is abelian $\forall k\in\Ze$
\item{$(vii)$} 
\quad $A_{k-1,k+2}=B_k+B_{k+1}$ is a factor $\forall k\in\Ze$
\end{Prop}

\begin{proof} By Lemma \ref{(1.4)} and eq. (\ref{equ(1.1)}), 
any vector in $M_{k+p}$ can be (uni\-quely) decomposed in sum of a
vector in $M_{k+p-1}$ and a vector in $B_{k+j-1}$. Iterating this
argument we get $(i)$.

Now we note that, by the global $j_k$-invariance of $B_k$,
$$
B_k\subseteq M_{k+1}\cap j_kM_{k+1}=M_{k+1}\cap M'_{k-1}=A_{k-1,k+1}.
$$ 
In particular, $B_{k+j}$ commutes with
$M_{k-1}$ for each $j\geq0$. Therefore, applying the decomposition
$(i)$ to a vector in $A_{k-1,k+p}$, we get $(ii)$. 

Let us take $x,y\in B_k$. Since $B_k$ is $j_k-$antiinvariant,
$$
Im(x,y)=Im(x,-j_ky)=Im(y,-j_kx)=Im(y,x)=-Im(x,y),
$$
which implies $(iii)$.

Now observe that
\begin{align*}
B_{k+1}\cap B'_k&=(B_{k+1}\cap M'_k)\cap B'_k=B_{k+1}\cap(M_k\vee
B_k)'=\cr &=B_{k+1}\cap M'_{k+1}=(M_{k+1}\cap M'_{k+1})\cap
Ker(j_{k+1}+I)\cr
\end{align*}
Since the elements in the center of $M_{k+1}$ are $j_{k+1}$
invariant  (see Remark 1.7, \cite{FiGu2}), it turns out that
$B_{k+1}\cap B'_k={0}$.

Then notice that $x\in\cz(\sum_{j=1}^{p}B_{k+j})$ {\it iff}
$x\in\sum_{j=1}^{p}B_{k+j}$ and
\begin{equation}\label{center}
Im(x,y)=0\ \forall y\in B_{k+j},\qquad1\leq j\leq p.
\end{equation}
Since any $x\in\sum_{j=1}^{p}B_{k+j}$ may be uniquely written as
$x=\sum_{1}^{p}b_{k+j}$, $b_{k+j}\in B_{k+j}$, and $B_m$ commutes with
$B_n$ if $|m-n|\geq2$, eq. (\ref{center}) for $j=1$ implies
$b_{k+2}=0$. Then eq. \ref{center} for $j=3$ implies $b_{k+4}=0$
and, iterating, $b_{k+2j}=0$ for any $1\leq j\leq p/2$. Analogously we
show that $b_{k+p-1}$ vanishes, and therefore $b_{k+p-2j-1}=0$ for
$j\leq p/2-1$. If $p$ is even this shows that $x=0$, i.e. $(iv)$.

When $p=2m+1$ is odd, we proved that $x=\sum_{0}^{m}b_{k+2j+1}$. Now,
taking $y\in B_{k+2j}$, $1\leq j\leq m$, and making use of the
$j$-anti-invariance of $B$, we have
\begin{align*}
0 &=Im(x,y)=Im(b_{k+2j-1}+b_{k+2j+1},y)\cr
  &=Im(b_{k+2j-1},y)-Im(b_{k+2j+1},j_{k+2j}y)\cr
  &=Im(b_{k+2j-1},y)+Im(j_{k+2j}b_{k+2j+1},y)\cr
  &=Im(b_{k+2j-1}+j_{k+2j}b_{k+2j+1},y),
\end{align*}
hence $b_{k+2j-1}+j_{k+2j}b_{k+2j+1}=0$, and this implies $(v)$.

Finally we observe that the inclusion $M_0\subset M_1$ is
irreducible \sse $A_{0,1}$, and therefore $A_{k,k+1},\quad\forall
k\in\Ze$, are trivial. Then 
$(vi)$ follows from $(ii)$ and $(iii)$ and
$(vii)$ follows from $(ii)$ and $(iv)$.
\end{proof}

\section{Inclusions of second quantization algebras}

We first recall that an inclusion of von~Neumann algebras  
$\cam_0\subset \cam_1$ is said of ${\rm depth}\,2$ if the von Neumann
algebra $\cam_0'\cap \cam_3$ is a factor.

\begin{Thm}\label{(2.1)} Let $M_0\subset M_1$ be an inclusion of
standard subspaces.
\item{$(i)$} If the codimension of $M_0$ in $M_1$ is finite and equal to $n$
$$
\R(M_1) = \R(M_0)\times_\alpha \br^n
$$
where $\br^n$ is identified with $B_0$ and $\alpha_h(\cdot)={\ad W(h)}$, $h\in
B_0$. 
\item{$(ii)$} If $M_0\subset M_1$ is irreducible then 
$\R(M_0)\subset\R(M_1)$ is of depth two.
\end{Thm}

\begin{proof} $(i)$ First we note that, by eq. \ref{CCR},
\begin{equation}\label{act}
\ad W(h)\left(W(k)\right)=e^{-iIm(h,k)}W(h+k)
\end{equation}
i.e. each closed  commutative subspace of $\H$ considered as an
additive group acts on any second quantization
algebra via $\ad W(\cdot)$. 

Then we only have to check that the hypotheses of Landstad
Theorem \cite{Land1} are satisfied. By Proposition \ref{(1.1)} and
Proposition \ref{(1.5)} ($i$) we get
$$
\R(M_1)=\R(M_0)\vee\{W(h):h\in B_0\}''
$$
i.e. $\R(M_1)$ is generated by $\R(M_0)$ and a unitary
representation of $B_0$ which is strongly continuous by the
finite-dimensionality of $B_0$. 

Now we claim that the action
$$
\b:h\in B_1\to \ad W(h)\in Aut(R(M_1))
$$
is dual to the action $\alpha$ of $B_0$ on $M_0$.

On the one hand
$$
\R(M_1)^\b=\R(M_1)\cap\{W(h):h\in B_1\}'=\R(M_1\cap B_1')=\R(M_0)
$$
where the last equality follows from Proposition \ref{(1.5)} $(i)$,
$(iv)$.  On the other hand, by eq. (\ref{act}), $B_1$ acts as the
dual group on $\R(B_0)$ \sse the pairing
$$
Im(\cdot,\cdot):B_0\times B_1\to\Re
$$
is a duality between real Hilbert spaces.
The factoriality of $B_0\vee B_1$ implies that the pairing is non
degenerate, the finite-dimensionality of $B_0$ and $B_1$ implies
that it is continuous, and the thesis follows.

$(ii)$ follows immediately by Proposition \ref{(1.5)}
\end{proof}

\subsection{Non regular inclusions}\label{4.1}

\begin{Thm}\label{(e1)} There exists an irreducible inclusion 
$M_0\subset M_1$ of standard subspaces with infinite codimension for
which the third relative commutant is of type III.
\end{Thm}

\begin{proof}
Let $M_0^n\subset M_1^n$ be a sequence of irreducible inclusions of
standard subspaces such that the codimension of $M_0^n$ in $M_1^n$ is
one, and let $M_0$, resp. $M_1$ the direct sum of the $M_0^n$,
resp. $M_1^n$ in $\oplus_{n=1}^\infty\H$.

As it is well known, direct sums of complex orthogonal real subspaces
give rise to tensor products at the second quantization level,
therefore the decomposition $A_{0,3}=\oplus_{n=1}^\infty A_{0,3}^n$
gives rise to an ITPFI decomposition. Since the codimension of $M_0^n$
in $M_1^n$ is one, $A_{0,3}^n$ is two-dimensional and the angle of
$A_{0,3}^n$ with $iA_{0,3}^n$ is a number $\theta_n$.

Therefore, by \cite{ArWo1} and \cite{FiGu2} Proposition 2.7, the type
of $A_{0,3}$ depends only on the sequence $\theta_n$. Choosing
isomorphic inclusions $M_0^n\subset M_1^n$ one gets a constant sequence
$\theta_n$ hence the factor has type III$_\lambda$, with
$\lambda=\tan^2\theta/2$.
\end{proof}

Let us recall that an inclusion $\cam_0\subset\cam_1$ is called {\it
regular} if the space of intertwiners $T:$ $yJ_1xJ_1T=TyJ_2xJ_2$ for
any $x\in\cam_1$, $y\in\cam_0$ is non zero. It was shown in
\cite{EnNe1} that for any depth two regular inclusion the third
relative commutant is of type I, and that for any crossed--product
inclusion regularity holds. This shows that the inclusions of Theorem
\ref{(e1)} are not regular, hence cannot come from a crossed--product
with a locally compact group. Indeed they constitute an example of a
kind of crossed--product with the group given by an infinite
dimensional Hilbert space with the additive structure.

We conclude this paper observing that, for the described inclusions,
the only visible invariant is the type of the third relative
commutant.  Since these inclusions are, in a sense, topologically
trivial, it seems reasonable to conjecture that the third relative
comutant is actually a complete invariant.

\begin{ack}
The first draft of this work was completed when the authors where
visiting the mathematics department of the Penn State University. We
thank the people in the department for the warm hospitality and Adrian
Ocneanu for having suggested this line of research.
\end{ack}

 
\end{document}